\documentclass[11pt]{amsart}

\usepackage{amssymb}

\usepackage{amsmath,latexsym,amsthm}

\usepackage{amsfonts}
\usepackage[greek,english]{babel}

\newtheorem*{theorem*}{Theorem}
\newtheorem*{corollary*}{Corollary}

\usepackage{hyperref}

\def\zetab{\bar\zeta}
\def\bel{\begin{equation}\label}
\def\eeq{\end{equation}}

\newtheorem{remark}{Remark}[section]
\newtheorem{definition}{Definition}[section]

\newtheorem{theorem}{Theorem}[section]
\newtheorem{lemma}[theorem]{Lemma}

\newtheorem{example}[theorem]{Example}

\newcommand\epf{\hfill{$\square$}\medskip}

\def\ds{\displaystyle}

\def\bega{\begin{array}}
\def\enda{\end{array}}

\def\bepmatrix{\begin{pmatrix}}
\def\enpmatrix{\end{pmatrix}}

\def\bel{\begin{equation}\label}
\def\eeq{\end{equation}}
\newcommand\ee{\end{equation}}
\def\benl{\begin{equation*}}
\def\eenl{\end{equation*}}

\def\forall{\hbox{for all }~}
\def\be{\begin{equation}}
\def\beq{\begin{equation}}
\def\bel{\begin{equation}\label}
\def\eeq{\end{equation}}

\newcommand\ba{\begin{array}}
\newcommand\ea{\end{array}}

\def\begi{\begin{itemize}}
\def\endi{\end{itemize}}

\def\pr{\partial}
\def\d{ {\rm d} }

\def\forall{\hbox{for all}~}

\newcommand{\cR}{\mathbb{R}}

\newcommand{\cN}{\mathbb{N}}

\def\C{\mathcal{C}}

\def\L{\mathcal{L}}

\def\xb{\bar{x}}

\def\zb{\bar{z}}

\def\xib{\bar{\xi}}
\def\etab{\bar{\eta}}
\def\zetab{\bar{\zeta}}

\def\uh{\hat{u}}

\def\eps{\varepsilon}

\title[Impulsive commutative systems]{A note on systems with ordinary and impulsive controls}

\begin{document}

\author[M.S. Aronna]{M. Soledad Aronna}
\address{M.S. Aronna, IMPA, Estrada Dona Castorina 110, Rio de Janeiro 22460-320, Brazil}
\email{aronna@impa.br}

\author[F. Rampazzo]{Franco Rampazzo}
\address{F. Rampazzo, Dipartimento di Matematica ,
Universit\`a di Padova\\ Padova  35121, Italy}
\email{rampazzo@math.unipd.it}

\thanks{This article was published in IMA J. Math. Control Inform. 2014.{\texttt{ doi:10.1093/imamci/dnu033}} \\
This work was partially supported by the European Union under the 7th Framework Programme FP7-PEOPLE-2010-ITN -  Grant agreement number 264735-SADCO, and the Fondazione CaRiPaRo Project
``Nonlinear Partial Differential Equations: models, analysis, and
control-theoretic problems".}

\maketitle

\begin{abstract}
{
We investigate an everywhere defined  notion of solution for  control systems whose dynamics  depend nonlinearly on the control $u$ and state $x,$ and are   affine in  the time derivative $\dot u.$ For this reason, the input $u,$ which is allowed to be Lebesgue integrable, is called {\it impulsive}, while a second, bounded measurable  control $v$ is denominated  {\it ordinary.} The proposed notion of solution  is derived from a topological (non-metric) characterization of a former concept of solution which was given  in the case when the drift  is $v$-independent. Existence, uniqueness and representation of the solution are studied,  and  a close analysis of effects of (possibly infinitely many) discontinuities on a null set is performed as well.
}

\vspace{6pt}

\end{abstract}

\section{Introduction}\label{SecIntro}

 Control systems of the form
\begin{align}
 \label{E}\tag{E}&\dot x = {f}(t,x,u,v) +\sum_{\alpha=1}^m g_\alpha(x) \dot u_\alpha,\quad \text{on } [a,b], \\
\label{IC}\tag{IC} &x(a)=\xb,
\end{align}
can be given  a classical interpretation as soon as  the control $u$ is an absolutely continuous function and the control $v$ is Lebesgue integrable.
This paper is devoted to the investigation of a notion of solution for the Cauchy problem (E)(IC), when one  assumes the following hypotheses:
\begin{itemize}
 \item[(i)] the vector fields $g_\alpha$ {\it commute}, namely $[g_\alpha,g_\beta]\equiv 0,$ for all $\alpha,\beta=1,\dots,m,$ where $[\cdot,\cdot]$ denotes the Lie bracket;
\item[(ii)] the inputs $u$ belong to  the space  $ \L^1([a,b];U)$ of  everywhere defined Lebesgue integrable functions.
\end{itemize}
Loosely  speaking, the denomination ``impulsive'' comes from the fact that, due to the affine dependence of the dynamics on the control's derivative $\dot u,$ a discontinuity in $u$ may cause a discontinuity in the corresponding trajectory $x.$
  On the other hand, the bounded, measurable  input  $v$ can be regarded as an ``ordinary" control.


Let us observe that the case where $u$ is taken in the class of bounded variation functions (and the commutativity  in (i) is not necessarily verified)  has received most of the attention (see  e.g. \cite{Ris65}, \cite{BreRam88}, \cite{DalRam91}, \cite{SilVin96} and references therein). 
In these articles, the authors studied the technique that is nowadays known as {\em graph completion.} An extension of this concept, also dealing with trajectories  with bounded variation, was investigated in  \cite{Kar06}, \cite{AruKarPer11} for systems of the form \eqref{E}, while a more general framework allowing the dependence of the vector fields $g_\alpha$ on the ordinary control $v$, was analyzed in \cite{AruKarPer12}.

 Even in the case where $u$ can have unbounded variation, a  notion of solution valid for systems where  $f$ is independent of the ordinary control $v,$   and both (i) and (ii) are met, has already been investigated  (see e.g. \cite{BreRam91}, \cite{Sar91}, \cite{Dyk94}). This solution can be  defined pointwise and verifies nice properties of uniqueness and continuity on the data.
The main goal of the present  note consists in investigating a  suitable generalization of this concept of solution to the case when $f$ is actually $v$-dependent.
 Incidentally, let us observe that  a system like
\bel{Eu}
\dot{x} = {f}(t,x,u,v) +\sum_{\alpha=1}^m g_\alpha(x,u) \dot u_\alpha,
\eeq
reduces to \eqref{E} as soon as one adds $m$ extra state variables $z_{1},\dots,z_{m}$ and the additional equations
$$
\dot{z}_{\alpha} = \dot u_\alpha,\quad \text{for }\alpha=1,\dots,m.
$$
In this case the commutative hypothesis (i) reads:
(i')  $\left[ g_\alpha+\frac{\partial}{\partial z^\alpha},g_\beta+\frac{\partial}{\partial z^\beta}\right]\equiv 0,$ for all $\alpha,\beta=1,\dots,0,$ where $[\cdot,\cdot]$ denotes the Lie bracket for vector fields on $\cR^{n+m}.$

In fact, several applications justify the introduction in the dynamical equations of the ordinary, bounded,  control $v$ besides the {impulsive} control $(u,\dot u).$ For instance, in Lagrangian mechanics, if the control  $u$ denotes the shape of a concatenation $C$  of rigid bodies and the input $v$ is, say, an external force or torque acting on  $ C,$ then the whole motion of  $C$ in space is determined by equations of the form \eqref{Eu}.  More generally, in a $N+m$-dimensional  Lagrangian system (where $N=n/2$)  the input $u$ might represent a portion of a local system of coordinates $(q,u),$  while  $x$ would be identified with $(q,p),$  $p$ being the {\it momenta} corresponding to the free coordinates $q$  (see \cite{BreAldo89}, \cite{Ram99}).
Let us point out that the commutativity assumption is actually verified  is some situations of practical interest \cite{BhaTiw09}.

The main results of the paper, including   existence, uniqueness, continuous dependence of solutions on data, state-response measure-zero changes of $u$, are stated in Section \ref{SecLimit}. The latter is concluded by Theorem \ref{TRep}, where a  representation of solutions is given  in terms of a diffeomorphism constructed through an application of  the Multiple Flow-box Theorem to the vector fields $\{g_1,\dots,g_m\}.$ All proofs can be found in Section  \ref{SecProofs}.

\vspace{3pt}

\noindent\textbf{Notation and  assumptions.}
 Let $h$ be a locally Lipschitz vector field on $\cR^n,$ and let $\xb \in \cR^n.$ Whenever the solution to
\benl
\dot x(t) = h(x(t)),\quad h(0)=\xb,
\eenl
is defined on an interval $I$ containing $0,$  we use
${\rm exp}({t h})(\xb)$ to denote the value of this solution at time $t.$

Let $I$ be a closed interval and let  $E$ be a subset  of an Euclidean  space $\cR^d.$ We use $\L^1(I;E)$ to denote the set of pointwise defined Lebesgue integrable functions from $I$ to $\cR^d$ with values in $E,$ while $L^1(I;E)$  will denote the corresponding  family of equivalence classes (with respect to the Lebesgue measure). We write $AC(I;E)$ for the set of absolutely continuous maps from $I$ to $E.$
For an open subset $\Omega \subseteq \cR^n,$ $\C^k(\Omega;\cR^d)$ will denote  the space of $k-$times continuously differentiable  $\cR^d$-valued functions defined on $\Omega.$

\vspace{3pt}

Throughout the paper we shall assume the following hypotheses on the control system \eqref{E}-\eqref{IC}:

\vspace{3pt}

\noindent {\bf Hypothesis  H:}
\begin{itemize}
\item[(i)]
$U$ is a compact subset of $\cR^m$ such that, for every bounded interval $I \subset \cR,$ for each $\tau \in I,$ and for every  function $u\in\L^1(I;U),$ there exists a sequence $(u_k^\tau) \subset AC(I;U)$ verifying
  $$
|u^\tau_k(\tau)- u(\tau)| +  \|u^\tau_k-u\|_1\to 0,
  $$
  when $k \to \infty.$\\
 \item[(ii)]
The set $V \subset \cR^l$ is compact.
 \item[(iii)] The map $f:[a,b]\times \cR^n \times \cR^m \times V \to \cR^m$
 is such that,
 \begin{itemize}
 \item[-] for each $(x,u,v)\in  \cR^n \times \cR^m \times V,$ the map $t\mapsto f(t,x,u,v)$ is measurable on $[a,b];$
 \item[-] for each $t\in [a,b],$ the function $(x,u,v) \to f(t,x,u,v)$ is continuous on $ \cR^n \times \cR^m \times V$ and, moreover,
 \item[-]  the map
 $$
 (x,u) \mapsto f(t,x,u,v),
 $$
 is locally Lipschitz on $\cR^n \times \cR^m,$ uniformly in $(t,v) \in [a,b] \times V.$  \end{itemize}
 \item[(iv)] For every $\alpha=1,\dots,m,$ $g_\alpha\in\C^1(\cR^n;\cR^n).$
  \item[(v)] There exists $A>0$ such that
  $$
  \left| \Big(f(t,x,u,v),g_1(x),\dots,g_m(x)\Big)\right| \leq  A(1+|(x,u)|),$$ for every $(x,u) \in \cR^n\times \cR^m$ uniformly in $(t,v) \in [a,b]\times V.$
  \end{itemize}

Notice that  hypotheses (ii)-(v) above imply that, for every initial value $\xb \in \cR^n,$ and each pair $(u,v)\in AC([a,b];\cR^m) \times L^1([a,b];V),$ the Cauchy problem \eqref{E}-\eqref{IC} has a unique (Carath\'eodory) solution, here denoted by $x[\xb,u,v].$
\vskip0.3truecm
\noindent {\bf Hypothesis  CC:}
\begin{itemize}
\item[(CC$_1$)] the vector fields $g_\alpha$ are {complete}\footnote{We say that $g_\alpha$ is {\em complete} if the solution to the Cauchy problem $\dot x=g_\alpha(x),\ x(0)=\xb \in \cR^n$ is (uniquely) defined on $\cR.$} and
\item[(CC$_2$)] $ g_1,\dots g_m$ verify the {\em global commutativity hypothesis} on $\cR^n,$
namely  for every Lipschitz continuous loop
$$
u:[0,1]\to \cR^m, \qquad  u(0)=u(1),
$$  and  each $\bar x\in \cR^n$ such that there exists a (unique)  Carath\'eodory solution to the Cauchy problem
$$
\dot x(t) = \sum_{\alpha=1}^m g_\alpha(x(t)) \dot u_\alpha(t),\quad t\in [0,1], \qquad x(0) =\bar x; $$
the solution $x$ is a loop, that is, it verifies  $x(0)=x(1)=\bar x.$
\end{itemize}

\begin{remark}{\rm
Let us define the {\em Lie bracket} of $g_\alpha$ and $g_\beta$ as
$$
[{g}_\alpha,{g}_\beta] := \ds\sum_{i=1}^n  \sum_{j=1}^n \left( \frac{\partial g_{\beta,i}}{\partial x_j} g_{\alpha,j} -  \frac{\partial g_{\alpha,i}}{\partial x_j} g_{\beta,j} \right)\frac{\partial}{\partial x_i}.
$$
It is trivial to verify that  for the domain $\cR^n$  the {\it null bracket condition}
\bel{lie0}
[{g}_\alpha,{g}_\beta] \equiv 0, \quad \alpha,\beta=1,\dots,m,
\eeq
 is necessary and sufficient for $g_1,\dots,g_m$ to verify the global commutativity hypothesis.
Actually, if instead of $\cR^n$ one considered an open subset $\Omega \subset \cR^n$ (or a differential manifold) as state space, the null bracket condition \eqref{lie0} would be no longer sufficient for global commutativity. As a trivial example, one can take
the vector fields  $g_1:= \left(1,0,\frac{-x_2}{x_1^2+x_2^2} \right)^\top,$ $g_2:=\left( 0,1,\frac{x_1}{x_1^2+x_2^2} \right)^\top,$ which verify the null bracket condition \eqref{lie0} on $\Omega:=\cR^2 \backslash \{0\},$ but  {\it do not} match the global commutativity hypothesis. Indeed, if
$
u(t):= \big(\cos (2\pi t),\sin (2\pi t) \big)^\top,$ for $t\in [0,1],$
 and $x$ is the  corresponding solution to
 $
 \dot x=g_1(x)\dot u_1 + g_2(x)\dot u_2,$ $x(0)=(1,0,0)^\top,$
one has
$x(1) = (1,0,2\pi)^\top\neq x(0).$}
\end{remark}



\section{Limit solutions}\label{SecLimit}

In this section we give the definition of limit solution and state the main results. The corresponding proofs have been  placed in  Section \ref{SecProofs}.

\begin{definition}[Limit Solution]\label{edsdef}
Consider an initial data $\xb \in \cR^n$ and controls  $(u,v)\in\L^1([a,b];U)\times L^1([a,b];V).$
We say that an $\L^1-$map $x:[a,b]\to \cR^n$ is a {\em limit solution} of the Cauchy problem \eqref{E}-\eqref{IC} if, for every $\tau\in [a,b],$ there exists a sequence $(u^\tau_k) \subset AC([a,b];U)$ such that:
\be
\label{limls}
|(x^\tau_k,u^\tau_k)(\tau)- (x,u)(\tau)|+ \|(x^\tau_k,u^\tau_k)- (x,u)\|_1\to 0,
\ee
where $x^\tau_k:=x[\xb,u^\tau_k,v].$
\end{definition}

\begin{remark}
Let us point out that $x$ is a limit solution associated to $u$ if, for every $\tau\in [a,b]$,   $(x,u)$ can be approximated, in the sense of \eqref{limls}, by sequences of absolutely continuous paths $(x_k^\tau,u_k^\tau)$   that verify \eqref{E} in the classical, Carath\'eodory sense.
We also observe that  no direct distributional meaning can be given to the derivative of $u$ or to \eqref{E}  (see some general considerations on the subject in  \cite{Haj85}), essentially because of two facts: on one hand the $g_\alpha$ are not constant; on the other hand we look for everywhere defined solutions.
\end{remark}

\begin{theorem}[Existence and uniqueness]
\label{Existence}
For every $\xb \in \cR^n,$ and every control pair $(u,v)\in\L^1([a,b];U)\times L^1([a,b];V),$ there exists a unique limit solution of the Cauchy problem \eqref{E}-\eqref{IC} defined on $[a,b].$
\end{theorem}

Given $\xb \in \cR^n,$ and a control pair $(u,v)\in\L^1([a,b];U)\times L^1([a,b];V),$ let $x[\xb,u,v]$ denote the (unique) corresponding limit solution of \eqref{E}-\eqref{IC}.

\begin{remark}
{\rm
Notice that,  for every input  $u\in\L^1([a,b];U),$ the map $t\mapsto \xb + u(t)-u(a)$ is a limit solution of the trivial Cauchy problem
\beq\label{trivial}
\dot x = \dot u,\quad   x(a) = \xb,
\eeq
and, thanks to the above uniqueness result, it is in fact  the only solution.
Since in general  $\L^1$ functions cannot be pointwise approximated by absolutely continuous functions (see e.g. to \cite{Ox80}),  the fact that the choice of the approximating control sequence depends on the time $\tau$ is crucial for guaranteeing existence of everywhere defined solutions, even for the  trivial equation \eqref{trivial}.
 }
\end{remark}

\vspace{5pt}

\begin{example}\label{minimum}{\rm

Let  $R\subset\cR^2$ be the subset defined by
$$
R\doteq \Big\{ (x,\hat y(x))\quad x\in[0,1]\Big\}\cup \Big\{(x, e^{-\frac{1}{2}})\quad x\in [2,3]\Big\}
$$
where
\benl
\hat y(x):=
\left\{
\ba{cl}
 e^x, & \text{for } x\in[0,\frac{1}{2}[,\\
 e^{1/2} e^{-2}, & \text{for } x \in \bigcup_{k=1}^{\infty} [1-\frac{1}{2k},1-\frac{1}{2k+1}[\,,\quad k\in\cN,\\
 e^{1/2} & \text{for } x \in \bigcup_{k=1}^{\infty} [1-\frac{1}{2k+1},1-\frac{1}{2k+2}[\,,\quad k\in\cN,\\
 e^{-1/2} & \text{for } x=1,
\ea
\right.
\eenl
and let us consider the  optimal  control problem\bel{opt}
\min \left\{w(2) + \big(y(1)- e^{-1/2}\big)^2 + \big(x(2)-3\big)^2\right\},
\eeq on the interval $[0,2]$ subject to  the dynamics
\bel{eqex1}
\left\{\begin{array}{l}
\dot x = 1 + \dot u_2, \\
\dot{y} = yv+ y\dot u_1,\\
\dot w =  \d\Big((x,y),R\Big),\\
(x,y,w)(0)= (0,1,0),
\end{array}\right.\eeq
 where the $v\in \{0,1\}$ and $(u_1,u_2)\in [-1,1]\times [0,1]$. Notice that \eqref{eqex1} meets the general  hypotheses, for the vector fields
$$
g_1\doteq \left(\begin{array}{c}1\\y\\0\end{array}\right), \quad g_2\doteq \left(\begin{array}{c}1\\0\\0\end{array}\right), \quad f\doteq \left(\begin{array}{c}1\\yv\\ \d\Big((x,y),R\Big)\end{array}\right),
$$
are Lipschitz continuous, and, moreover,  $g_1,g_2$ are smooth and verify $[g_1,g_2]\equiv 0$.

We claim that the limit solution $(x,y,w)$ corresponding to the input
\benl
v(t):=
\left\{
\ba{cl}
1 & \text{for } t\in[0,1/2[,\\
0 & \text{for } t\in [1/2,1],
\ea
\right.
\eenl
\benl
u_1(t):=
\left\{
\ba{cl}
(-1)^{k+1} & \text{for } t\in[1-\frac{1}{k},1-\frac{1}{k+1}[,\,\, k \in\cN,\\
0 & \text{for } t\in [1,2],
\ea
\right.
\eenl
\benl
u_2(t):=
\left\{
\ba{cl}
0 & \text{for } t\in[0,1],\\
1 & \text{for } t\in ]1,2],
\ea
\right.
\eenl
is a minimum for problem \eqref{opt}. Indeed, on any subinterval $[t_1,t_2]\subset[1/2,1]$ where $u_1$ is absolutely continuous, one has
\be
\label{forex}
y(t)=y(t_1)e^{u_1(t) - u_1(t_1)}.
\ee
 On the other hand, one can easily check that
\benl
\begin{split}
y(1-1/k  +) &= y(1-1/k -)e^2, \quad \text{if $k$ is odd},\\
y(1-1/k  +) &=  y(1-1/k -)e^{-2}, \quad \text{if $k$ is even,}
\\
y(1) &=e^{-1/2},
\end{split}
\eenl
where $y(1-1/k -)$ and $y(1-1/k  +)$ denote the left and right limits of $y$ at $t=1-1/k,$ respectively. Moreover, $$
x= t\,\,\,\forall t\in [0,1], \quad x=t+1\,\,\, \forall t\in ]1,2].$$  Hence  $\d\Big((x(t),y(t)),R\Big)=0$ for all $t\in[0,2]$  and $x(2)=3$, so the corresponding  payoff is equal to zero. Therefore $(x,y,w)$ is an optimal trajectory, since the payoff of every control-trajectory pair is nonnegative.

Notice that both $(u_1,u_2) $ and $(x,y,w)$ have infinitely many discontinuities and  unbounded variation. Observe also that it is crucial that the input $u$ and the  solution are defined everywhere. In fact, the control $(v,\tilde u_1, u_2)$, with $\tilde u_1(t) = u_1(t)$ for all $t\neq 1$ and $\tilde u_1(1)=1$, is not optimal, for the corresponding solution $(\tilde x,\tilde y,\tilde w)$ is equal to $(x,y,w)$ on  $ [0,2]\backslash \{1\},$ while, in view of Theorem \ref{negl-jumps}, $\tilde y(1) = e^{1/2}.$
}
\end{example}

\begin{theorem}[Continuous dependence]
 \label{cont-depL1}
 The following assertions hold true:
\begin{itemize}
\item[{(i)}] for each $\xb \in \cR^n$ and $u\in \L^1 ([a,b];U),$ the function $v\mapsto x[\xb,u,v]$ is continuous from $L^1([a,b];V)$ to $L^\infty([a,b];\cR^n);$
\item[(ii)] for any $r>0,$ there exists a compact subset $K'\subset \cR^n,$ such that the trajectories $x[\bar x,u,v]$ have values in $K',$ whenever we consider $|\xb|\leq r,$  $u\in \L^1([a,b];U)$ and $v\in L^1([a,b];V) ;$
\item[(iii)] for each $r>0,$  there exists a constant $M>0$ such that, for every $\tau \in[a,b],$ for all $|\xb_1|, |\xb_2|\leq r,$ $u_1,u_2 \in \L^1([a,b];U)$  and for every $v \in L^1([a,b];V),$ one has
\bel{EstL1}
|x_1(\tau)-x_2(\tau)| +  \|x_1-x_2\|_1 \leq  
\, M\Big[ |\xb_1-\xb_2|+ |u_1(a) - u_2(a)|  + |u_1(\tau)-u_2(\tau)|+\|u_1-u_2\|_1
\Big],
\ee
where  $x_1 := x[\xb_1,u_1,v],$ $x_2 := x[\xb_2,u_2,v].$
\end{itemize}
\end{theorem}

Since the limit  solution depends on the pointwise definition of $u,$ it is interesting to investigate the effects of a change of the $u$'s values on a measure-zero subset of $[a,b].$

\begin{theorem}[Pointwise dependence]\label{negl-jumps}
Let us consider  an interval $[a,b],$ an initial state $\bar x \in \cR^n,$ and an ordinary control $v \in L^1([a,b];V).$ Let $u,\hat u \in \L^1([a,b];U)$ be impulse controls  that coincide a.e. in $[a,b]$ and that verify $u(a)=\uh(a).$
  Then, setting $x:= x[\xb, u,v],$ $\hat x:= x[\xb, \hat u,v],$ one has
\bel{jump}
x(t)=  \exp \left({\sum_{\alpha=1}^m (u_\alpha(t)-\hat u_\alpha(t))  g_\alpha}\right) \big(\hat x(t)\big), \quad\forall t\in [a,b].
\ee
In particular,
$$ x(t)=\hat x(t)$$
for every $t\in [a,b]$ such that $u(t)=\hat u(t),$ that is, almost everywhere.
\end{theorem}

\vskip0.4truecm
In order to state the representation theorem below we need to introduce a change of coordinates induced by the  $g_\alpha$'s flows.

Let us extend $f,g_\alpha,$ for $\alpha =1,\dots,m$  to functions $\tilde f,\tilde g_\alpha$  with values in  $\cR^{n+m}$ by setting, for every $(t,x,z,v)\in [a,b]\times  \cR^{n+m} \times V,$
$$
\tilde f(t,x,z,v):=\sum_{j=1}^n {f}_j(t,x,z,v)\frac{\partial}{\partial x_j},\quad
\tilde g_\alpha(x,z):= \sum_{j=1}^n {g}_{\alpha,j} (x,z) \frac{\partial}{\partial x_j} +     \frac{\partial}{\partial z_\alpha}\,,
$$
where $\left(\frac{\partial}{\partial x_1},\dots,\frac{\partial}{\partial x_n}, \frac{\partial}{\partial z_1},\dots, \frac{\partial}{\partial z_m}\right)$ is the canonical basis of $\cR^{n+m}.$
\footnote{Notice that at this stage there is no more need to distinguish between
 equations \eqref{E} and \eqref{Eu}, for the $\tilde g_\alpha$ are vector fields on $\cR^{n+m},$ so it is irrelevant that their first $n$ components are or are not dependent on $u.$}
Let  $\mathrm{Pr}:\cR^{n} \times \cR^m\rightarrow \cR^n$ denote the {canonical projection} on the first factor, i.e.
$$
\mathrm{Pr}(x,z) := x,
$$
and let the function  $\varphi:\cR^n\times \cR^{m}\to \cR^n$ be defined by
$$
\varphi(x,z) :=
\mathrm{Pr}\circ \exp \left(- z_{m} \tilde g_{m}\right)\circ\dots \exp \left(- z_{1} \tilde g_{1}\right)
(x,z).
$$

Finally, let us consider the map $\phi:\cR^n\times \cR^{m}\rightarrow \cR^n\times \cR^{m}$ given by
$$
\phi(x,z) := (\varphi(x,z),z).
$$

\begin{lemma}
Assume that the vector fields $g_1,\dots,g_m$ belong to  $\C^r(\cR^n;\cR^n),$ with $r\geq 1.$ Then the  mapping $\phi$ is a $\C^r$-diffeomorphism of $\cR^{n+m}$ onto itself and, for every $(\xi,\zeta)\in \cR^{n+m},$ one has
$$
\phi^{-1}(\xi,\zeta) = (\varphi(\xi,-\zeta), \zeta).
$$
\end{lemma}

The $\C^r$-diffeomorphism  $\phi$  induces the $\C^{r-1}$-diffeomorphism $D\phi$ on the tangent bundle, where $D$ denotes differentiation. For each $\alpha=1,\hdots,m,$ $t\in [a,b],$ $(\xi,\zeta,v)\in \cR^{n+m}\times V,$ let us set
$$
\begin{array}{l}
\tilde F(t,\xi,\zeta,v):= D\phi (x,z) \, \tilde f(t,x,z,v),
\\
\tilde G_\alpha(\xi,\zeta):= D\phi (x,z) \, \tilde g_{\alpha}(x,z),
\end{array}
$$
where $(x,z):= \phi^{-1}(\xi,\zeta).$ As a direct consequence of the  {\em Simultaneous Flow-Box Theorem}  (see e.g.  \cite{Lang}), one obtains the following result  (see Lemma 2.1 in \cite{BreRam91} for a proof).

\begin{lemma}
\label{flowbox}
For every $i=1,\dots, n$ and $\alpha=1,\dots, m$ one has
$$
\tilde F =\sum_{i,j=1}^n\left(\frac{\pr \phi_i}{\pr x_j} {f}_j \right) \frac{\pr}{\pr x_i},\quad \tilde{G}_{\alpha}
= \frac{\pr }{\pr z_\alpha},
$$
where we have set $\phi=(\phi_1,\dots,\phi_{n+m}).$
\end{lemma}

Notice that the last $m$ components of $\tilde F$ are zero. More precisely, $\tilde F$ can be written in components as $\tilde F=\begin{pmatrix} {F} \\ 0 \end{pmatrix}$ with $F:[a,b]\times\cR^{n+m }\times V \to \cR^n.$
Consider the Cauchy problem
\begin{align}
\label{Et}
\dot{\xi} &= {F}(t,\xi,u,v),\quad \text{on } [a,b],\\
\label{ICt}
\xi(a)&=\xib.
\end{align}
For each $\xib \in \cR^n,$ $(u,v)\in \L^1([a,b];\cR^n) \times L^1([a,b];V),$ there exists a (unique)  Carath\'eodory solution of \eqref{Et}-\eqref{ICt}, which will be  here denoted by $\xi[\xib,u,v].$

Theorem \ref{TRep} below , which is trivial in the case $u\in AC,$ provides a representation of the solutions  of \eqref{E} (and of \eqref{Eu}) in terms of images of solutions of the simpler equation \eqref{Et} through the map $\varphi$ previously introduced.
\begin{theorem}[Representation of limit solutions]
\label{TRep}
For any $\xb\in \cR^n,$ $(u,v) \in \L^1([a,b];U)\times L^1([a,b];V),$ one has
$$
\xi[\xib,u,v](t) = \varphi\big(x[\xb,u,v](t),u(t) \big),\quad \text{for all } t\in [a,b],
$$
where we have set $\xib:=\varphi(\xb,u(a)).$
\end{theorem}

The proof of this theorem is given in the next section.


\section{Proofs of the results of Section \ref{SecLimit}}\label{SecProofs}

Since we are going to exploit the diffeomorphism $\phi:  \cR^{n+m} \to \cR^{n+m}$ it is convenient to embed \eqref{E}-\eqref{IC} in the $n+m$-dimensional Cauchy problem
\begin{align}
\label{AS}
\begin{pmatrix}
\dot{x}\\
\dot{z}
\end{pmatrix}
&= \tilde f(t,x,z,v)+   \sum_{\alpha=1}^m \tilde{g}_\alpha(x,z) \dot  u_\alpha,\\
\label{ASIC}
\begin{pmatrix}
{x}\\{z}
\end{pmatrix}(a) &= \begin{pmatrix} \xb \\ \zb \end{pmatrix}.
\end{align}
 In view of the considered hypotheses, when $u\in AC([a,b];U),$ for every $(\xb,\zb) \in \cR^{n+m}$ and $v \in L^1([a,b];V),$  there exists a unique solution to \eqref{AS} in the interval $[a,b].$ We let  $(x,z)[\xb,\zb,u,v]$ denote this solution.

We shall also consider  the  Cauchy problem
\begin{align}
\label{TS}
\begin{pmatrix}
\dot{\xi}\\
\dot{\zeta}
\end{pmatrix}
&= \tilde F(t,\xi,\zeta,v)+  \sum_{\alpha=1}^m \tilde G_\alpha  \dot u^\alpha,\\
\label{TSIC}
\begin{pmatrix}
{\xi}\\{\zeta}
\end{pmatrix}(a)&= \begin{pmatrix} \xib \\ \zetab \end{pmatrix}
.
\end{align}
Also for this problem, when $u\in AC([a,b];U),$ for every $(\xib,\zetab)\in  \cR^{n+m}$ there exists a unique solution to \eqref{TS}-\eqref{TSIC} in the interval $[a,b].$ We let  $(\xi,\zeta)[\xib,\zetab,u,v]$ denote this solution.
\begin{remark}\label{equivAC}{\rm When $u\in AC([a,b];U),$  the relation between the two systems is given by
$$
(\xi,\zeta)[\xib,\zetab,u,v](t)  = \phi\Big( (x,z)[\xb,\zb,u,v](t) \Big),
\quad \text{for all } t\in [a,b],
$$
where $(\xib,\zetab) := \phi(\xb,\zb).$}
\end{remark}
\begin{remark}{\rm
The crucial difference between the two latter systems relies on the fact that the vector fields  $\tilde G_\alpha=\frac{\partial}{\partial z_\alpha}$ are constant}.
\end{remark}

Theorem \ref{ContAC-simple} below   will be utilized to prove Theorem \ref{cont-depL1}, of which it is in fact a particular case.

\begin{theorem}
\label{ContAC-simple}
The following assertions hold.
\begin{itemize}
\item[(i)] For each $(\xib,\zetab)\in \cR^{n+m}$ and $u\in AC([a,b];U),$ the function $v\mapsto \xi[\xib,\zetab,u,v]$ is continuous from $L^1([a,b];V)$ to $L^\infty([a,b];\cR^n).$
\item[(ii)] Furthermore, for any $r>0,$  there exists a compact subset $K'\subset  \cR^{n+m}$ such that the trajectories $(\xi,\zeta)[\xib,\zetab,u,v]$ have values in $K',$ whenever we consider $|(\xib,\zetab)|\leq r,$ $u\in AC([a,b];U),$ and $v\in L^1([a,b];V).$
\item[(iii)] Finally, for each $r>0,$ there exists a constant $M>1$ such that, for every $\tau \in[a,b],$ for all $|(\xib_1,\zetab_1)|, |(\xib_2,\zetab_2)| \leq r,$ for all $u_1,u_2 \in AC([a,b];U),$ and for every $v \in L^1([a,b];V),$ one has
\bel{EstL1xi}
\begin{split}
|(\xi_1,\zeta_1)(\tau)-&(\xi_2,\zeta_2)(\tau)| +  \|(\xi_1,\zeta_1)-(\xi_2,\zeta_2)\|_1 \leq \\
&M\Big[ |(\xib_1,\zetab_1)-(\xib_2,\zetab_2)|+ |u_1(a)-u_2(a)| + |u_1(\tau)-u_2(\tau)|
+\|u_1-u_2\|_1
\Big].
\end{split}
\ee
where $(\xi_1,\zeta_1) := (\xi,\zeta)[\xib_1,\zetab_1,u_1,v_1]$ and $(\xi_2,\zeta_2) := (\xi,\zeta)[\xib_2,\zetab_2,u_2,v_2].$
\end{itemize}
\end{theorem}

To prove Theorem \ref{ContAC-simple}, we shall exploit the following fixed-point result on parameterized contraction mappings (see e.g.  Theorem A.1 in\cite{BrePicBook}).

\begin{lemma}
\label{BC}
Let $X$ be a Banach space, $\Lambda$ a metric space and $\chi:\Lambda\times X \to X$ be a continuous function such that
\be
\label{BCeq1}
\| \chi(\lambda,x) - \chi(\lambda,y) \| \leq L \|x-y\|,\quad \forall \lambda\in \Lambda,\ x,y\in X,
\ee
with $L<1.$ Then the following assertions hold.
\begin{itemize}
\item[(a)] For every $\lambda \in \Lambda,$ there exists a unique $x(\lambda)$ such that
$$
x(\lambda) = \chi(\lambda,x(\lambda)).
$$
\item[(b)] The map $\lambda\mapsto x(\lambda)$ is continuous, and one has
$$
\| x(\lambda)-x(\tilde \lambda) \| \leq \frac{1}{1-L} \| \chi(\lambda,x(\tilde\lambda)) - \chi(\tilde\lambda,x(\tilde\lambda)) \|.
$$
\end{itemize}
\end{lemma}

\vspace{5pt}


{\em Proof of Theorem \ref{ContAC-simple}.}
Item (i) follows from classical results of continuity of the input-output map of a control system.

 To prove the remaining assertions, assume momentarily  that $F$ is globally Lipschitz continuous with respect to the variable $(\xi,\zeta)$ with Lipschitz constant $L.$ Later we shall remove this extra assumption.

 For $(\xib,\zetab,u,v)\in \Lambda:=\cR^{n+m}\times AC ([a,b];U)\times L^1([a,b];V)$ and $(\xi,\zeta) \in X:= AC ([a,b];\cR^{n+m}),$ let us consider the mapping $\chi:\Lambda \times X \to X$ such that, for all $t\in [a,b],$
\benl
\chi(\xib,\zetab,u,v,\xi,\zeta)(t):= 
\begin{pmatrix} \xib \\ \zetab \end{pmatrix}
+ \int_a^t \tilde F(s,\xi(s),\zeta(s),v(s)) {\rm d} s
+\sum_{\alpha=1}^m [u_\alpha(t) - u_\alpha(a)]{\bf e}_{n+\alpha},
\eenl
where ${\bf e}_{n+\alpha}$ denotes the $(n+\alpha)$th vector of the canonical basis of $\cR^{n+m}.$
Observe that 
$$
\begin{pmatrix} \xi \\ \zeta \end{pmatrix}=\chi(\xib,\zetab,u,v,\xi,\zeta) \text{ if and only if  }
\begin{pmatrix} \xi \\ \zeta \end{pmatrix} \text{ is solution of \eqref{TS}.}
$$
 We are therefore interested in applying the fixed-point result in Lemma \ref{BC} to the function $\chi.$

Fix $\tau \in [a,b],$ and in the space $\Lambda$ consider the norm
$$
\|(\xib,\zetab,u,v)\|_Y := |(\xib,\zetab)|+|u(a)| + |u(\tau)| + \|u\|_1+\|v\|_1,
$$
and in $X,$  define the norm
\benl
\| (\xi,\zeta) \|_X :=
\frac{e^{-4(b-a)L}}{4L} |(\xi,\zeta)(\tau)|
+\int_a^b e^{-4sL}|(\xi,\zeta)(s)|{\rm d}s.
\eenl
We shall prove that $\chi$ is continuous from $(\Lambda\times X,\|\cdot\|_Y+\|\cdot\|_X)$ to $(X,\|\cdot\|_X).$
By the Lipschitz continuity of the maps $(\xi,\zeta) \to \tilde F(t,\xi,\zeta,v),$ for any $(\xib_1,\zetab_1,u_1,v_1,\xi_1,\zeta_1),$  $(\xib_2,\zetab_2,u_2,v_2,\xi_2,\zeta_2)$  in $\Lambda\times X,$ one has
\be
\label{estchi1}
\begin{split}
&\|\chi (\xib_1,\zetab_1,u_1,v_1,\xi_1,\zeta_1)-\chi(\xib_1,\zetab_1,u_1,v_1,\xi_2,\zeta_2)\|_X =\\
& \frac{e^{-4(b-a)L}}{4L} \left| \int_a^\tau \big[\tilde F \big(s,\xi_1(s),\zeta_1(s),v_1(s) \big) - \tilde F \big(s,\xi_2(s),\zeta_2(s),v_1(s) \big) \big] {\rm d}s \right| +\\
& \int_a^b e^{-4sL} \left| \int_a^t \big[\tilde F \big(s,\xi_1(s),\zeta_1(s),v_1(s) \big) - \tilde F \big(s,\xi_2(s),\zeta_2(s),v_1(s) \big) \big] {\rm d}s \right| {\rm d} t \leq \\
&\qquad \frac{1}{4} \| (\xi_1,\zeta_1) - (\xi_2,\zeta_2) \|_X,
\end{split}
\ee
\be
\label{estchi2}
\begin{split}
&\|\chi (\xib_1,\zetab_1,u_1,v_1,\xi_2,\zeta_2)-\chi(\xib_2,\zetab_2,u_2,v_2,\xi_2,\zeta_2)\|_X =\frac{e^{-4(b-a)L}}{4L} \left| \begin{pmatrix} \xib_1-\xib_2 \\ \zetab_1-\zetab_2 \end{pmatrix}  \right| + \\
 & \frac{e^{-4(b-a)L}}{4L} \left| \int_a^\tau \big[\tilde F \big(s,\xi_2(s),\zeta_2(s),v_1(s) \big) - \tilde F\big(s, \xi_2(s),\zeta_2(s),v_2(s) \big) \big] {\rm d}s \right|  +\\
& \int_a^b e^{-4sL} \left| \int_a^t \big[\tilde F \big(s,\xi_2(s),\zeta_2(s),v_1(s) \big) - \tilde F \big(s,\xi_2(s),\zeta_2(s),v_2(s) \big) \big] {\rm d}s \right| {\rm d} t +\\
&\frac{e^{-4(b-a)L}}{4L}\left| \sum_{\alpha=1}^m \big[ u_{1,\alpha}(\tau)-u_{1,\alpha}(a)-u_{2,\alpha}(\tau)+u_{2,\alpha}(a)\big] {\bf e}_{n+\alpha}   \right| +\\
&\int_a^b e^{-4sL} \left| \int_a^t  \sum_{\alpha=1}^m \big[ u_{1,\alpha}(s)-u_{1,\alpha}(a)-u_{2,\alpha} (s)+u_{2,\alpha} (a)\big] {\bf e}_{n+\alpha} {\rm d} s \right| {\rm d} t.
\end{split}
\ee
By the Dominated Convergence Theorem, for any fixed trajectory $( \xi, \zeta),$ the  mapping
$$
v \mapsto \tilde F(\cdot,\xi, \zeta, v )
$$
is continuous from $L^1([a,\tau];V)$ to $L^1([a,\tau]; \cR^{n+m}).$
Thus, by \eqref{estchi1}-\eqref{estchi2}, for each $(\xib_1,\zetab_1,u_1,v_1,\xi_1,\zeta_1)\in \Lambda \times X$ and for every $\eps>0,$ there exists $\delta>0$ such that, if
\benl
\begin{split}
&|(\xib_2,\zetab_2)-(\xib_1,\zetab_1)| +|u_2(a)-u_1(a)| \\
&\quad +|u_2(\tau)-u_1(\tau)|+\|u_2-u_1\|_1+ \|v_2-v_1\|_1 + \| (\xi_2,\zeta_2) - (\xi_1,\zeta_1) \|_X < \delta
\end{split}
\eenl
then
\benl
\|\chi (\xib_2,\zetab_2,u_2,v_2,\xi_2,\zeta_2)-\chi(\xib_1,\zetab_1,u_1,v_1,\xi_1,\zeta_1)\|_X <\eps.
\eenl
Hence, $\chi$ is continuous, and in  view of \eqref{estchi1} the inequality \eqref{BCeq1} holds true. To apply Lemma \ref{BC},   let us identify  $\lambda$ with $(\xib_1,\zetab_1,u_1,v_1)$ and $\tilde\lambda$ with $(\xib_2,\zetab_2,u_2,v_2).$ Then one has that there exist $(\xi_1,\zeta_1),(\xi_2,\zeta_2) \in X$ such that
$$
\chi(\xib_1,\zetab_1,u_1,v_1,\xi_1,\zeta_1) = (\xi_1,\zeta_1),\quad
\chi(\xib_2,\zetab_2,u_2,v_1,\xi_2,\zeta_2) = (\xi_2,\zeta_2),
$$
that is $(\xi_1,\zeta_1) = (\xi,\zeta)[\xib_1,\zetab_1,u_1,v_1]$ and $(\xi_2,\zeta_2) = (\xi,\zeta)[\xib_2,\zetab_2,u_2,v_2].$

In view of item (b) in Lemma \ref{BC}, we get
\bel{Contrac1}
 \left\|  \begin{pmatrix}  \tilde\xi - \xi \\ \tilde\zeta - \zeta  \end{pmatrix}  \right\|_X \leq \frac{1}{2L} \left[ |(\tilde \xib,\tilde \zetab)-(\xib,\zetab)| \right.  \left.\, + |\tilde u(a)-u(a)| + |\tilde u(\tau)-u(\tau)|\right] + \int_a^b e^{-4tL} |\tilde u(t)-u(t)| {\rm d}t.
\ee
Therefore, by the inequalities
\bel{Contrac2}
  \left\| \begin{pmatrix} \xi_2 - \xi_1 \\ \zeta_2 - \zeta_1  \end{pmatrix}  \right\|_X
 \geq \frac{e^{-4(b-a)L}}{4L} \left|  \begin{pmatrix} \xi_2(\tau) - \xi_1(\tau) \\ \zeta_2(\tau) - \zeta_1(\tau)  \end{pmatrix}  \right| + e^{-4bL}  \left\|  \begin{pmatrix} \xi_2 - \xi_1 \\ \zeta_2 - \zeta_1  \end{pmatrix}  \right\|_1,
\ee
and  \eqref{Contrac1}  one obtains \eqref{EstL1xi}, with a constant $M$ depending only on $L,$ and hence item (iii) is proved. Notice also that item (ii) is a consequence of the following standard result on ODE's:
\begin{lemma}[Bounds on solutions]
\label{Apriori}
Under the general hypothesis H, for each $r>0,$  there exists a compact set $K'\subset \cR^{n+m},$ such that any solution $(\xi,\zeta)$ of \eqref{TS}-\eqref{TSIC} remains in $K'$ whenever $|(\xib,\zetab)|\leq r,$ $u\in AC([a,b];U)$ and $v\in L^1([a,b];V).$
\end{lemma}
This completes the proof of the theorem  under the additional assumption that $F$ is globally Lipschitz.

\vskip0.2truecm

To prove the general case we use a standard cut-off function argument.

 Take $r>0,$ and let $K'\subset  \cR^{n+m}$ be the compact set provided by Lemma \ref{Apriori}. Let $\rho \in \C^1(\cR^{n+m})$ be a smooth real function such that $\rho=1$ on $K'$ and $\rho=0$ outside a neighborhood of $K'.$
Define $\hat F (t,\xi,\eta,v):= \rho(\xi,\eta) \tilde F(t,\xi,\eta,v),$ and set
$$
\Lambda:= \{(\xib,\zetab)\in \cR^{n+m}:|(\xib,\zetab)|\leq r\} \times AC([a,b];U) \times L^1([a,b];V).$$
Then, for $(\xib,\etab,u,v)\in \Lambda,$ the corresponding solution $(\xi,\eta)$ of the Cauchy problem
$$\begin{pmatrix}
\dot{\xi}\\
\dot{\zeta}
\end{pmatrix}
= \hat F(t,\xi,\zeta,v)+  \sum_{\alpha=1}^m \tilde G_\alpha  \dot u^\alpha,\quad
\begin{pmatrix}
{\xi}\\{\zeta}
\end{pmatrix}(a)= \begin{pmatrix} \xib \\ \zetab \end{pmatrix}
,
$$
coincides with $(\xi,\eta)[\xib,\etab,u,v]$ (and remains inside $K'$).
Now the function $\hat F$ is globally Lipschitz, and the procedure done before can be repeated for this new metric space $\Lambda$ and for the function $\hat F$ in the place of $\tilde F.$
Therefore, one  can obtain  the estimate \eqref{EstL1xi} with  a constant $M$ depending only on the Lipschitz constant of the mapping $(\xi,\eta) \mapsto \tilde F(t,\xi,\eta,v)$ in the set $K'.$ This completes the proof of Theorem \ref{ContAC-simple}.
\epf

{\em Proofs of Theorems \ref{Existence} and \ref{TRep}.}
Let $\xb\in \cR^n, $  $(u,v) \in \L^1([a,b];U)\times L^1([a,b];V),$ and let $\xi:=\xi[\xib,u,v]$ be  the unique Carath\'eodory solution of \eqref{Et} with the initial condition $\xi(a)=\xib := \varphi(\xb,u(a)).$  Define, on $[a,b],$ the function
\be
\label{xdef}
x:=\varphi\circ(\xi,-u),
\ee
and let us show that $x$ is the unique limit solution of \eqref{E}-\eqref{IC} associated with $\xb$ and $(u,v).$

Choose $\tau \in [a,b],$ and consider a sequence of absolutely continuous controls $u_k^\tau: [a,b] \to U$ verifying
\be
\label{uktau}
|u^\tau_k(a) - u(a)|+|u^\tau_k(\tau) - u(\tau)|+\|u_k^\tau-u\|_1\to 0.
\ee
Consider the equation \eqref{Et} with the initial condition
$$
\xi(a)=\xib_k^\tau:=\varphi(\xb,u^\tau_k(a)).
$$
 Let $\xi_k^\tau := \xi[\xib_k^\tau,u_k^\tau,v]$ be the (unique) corresponding Carath\'eodory solution. Then, by standard results of continuity with respect to the data, one has that
\be
\label{xikunif}
\xi_k^\tau \to \xi,\  \text{uniformly on }[a,b].
\ee
For the augmented system \eqref{TS}, $(\xi^\tau_k,u^\tau_k)$ is the unique solution with the initial conditions $\xi(a)= \xib_k^\tau,$ $\eta(a)=u^\tau_k(a).$ In view of item (iii) in Theorem \ref{ContAC-simple}, the functions $(\xi^\tau_k,u^\tau_k)$ have values in a compact set $\tilde{K}'\subset  \cR^{n+m}.$
In view of  Remark  \ref{equivAC}, the map
$$
(x_k^\tau,z_k^\tau):=\phi^{-1} \circ (\xi_k^\tau,u_k^\tau),
$$
is, for each $k\in \cN,$ the unique Carath\'eodory solution of \eqref{AS} with the initial conditions $x(a)=\xb,$ $z(a) = u_k^\tau(a).$ Notice that the functions $ (x_k^\tau,z_k^\tau)$ have values inside the compact set $K':=\phi^{-1}(\tilde K').$ In particular, the $x_k^\tau$'s are uniformly bounded.
Observe as well that
$$
|x(\tau)-x_k^\tau(\tau)| = |\varphi(\xi(\tau),-u(\tau))-\varphi(\xi_k^\tau(\tau),-u_k^\tau(\tau))| \to 0,
$$
due to  \eqref{xikunif}-\eqref{uktau} and the continuity of $\varphi.$  Furthermore, since $u_k^\tau \to u$ almost everywhere and all these functions are uniformly bounded, one gets
$$
\|x-x_k^\tau\|_1 = \int_a^b \big|\varphi(\xi(t),-u(t))-\varphi(\xi_k^\tau(t),-u_k^\tau(t)) \big| {\rm d} t
\to 0,
$$
thanks to the Dominated Convergence Theorem. Thus, $x$ is a limit solution of \eqref{E}-\eqref{IC} associated with $\xb$ and $(u,v).$

We shall now prove the uniqueness. Suppose on the contrary that there are two different limit solutions $x$ and $\tilde x$ of \eqref{E}-\eqref{IC} associated with $(\xb,u,v).$ Let $\tau \in ]a,b]$ be such that $x(\tau) \neq x^*(\tau),$ and let $(u_k^\tau),$ $({u}_k^{\tau,*})$ be sequences in $AC([a,b];U)$ verifying
\be
\label{uktau}
|u^\tau_k(a)- u(a)|+|(x^\tau_k,u^\tau_k)(\tau)- (x,u)(\tau)|+ \|(x^\tau_k,u^\tau_k)- (x,u)\|_1\to 0,
\ee
and
\be
\label{tildeuktau}
|u^{\tau,*}_k(a)- u(a)|+|( x^{\tau,*}_k,u^{\tau,*}_k)(\tau)- (x^{*},u)(\tau)|+ \|( x^{\tau,*}_k, u^{\tau,*}_k)- ( x^*,u)\|_1
\to 0,
\ee
where $x^\tau_k:=x[\xb,u^\tau_k,v],$ $ x^{\tau,*}_k:=x[\xb, u^{\tau,*}_k,v].$
Let $\xi_k^\tau := \xi[\xib_k^\tau,u_k^\tau,v],$ $\xi_k^{\tau,*} := \xi[\xib_k^{\tau,*},u_k^{\tau,*},v]$ be the solutions of \eqref{Et} with initial conditions
$$
\xi(a)=\xib_k^\tau:=\varphi(\xb,u_k^\tau(a)),\quad \xi(a)=\xib_k^{\tau,*}:=\varphi(\xb, u_k^{\tau,*}(a)),
$$
respectively.
Then, since both $(u_k^\tau),$ $({u}_k^{\tau,*})$ converge to $u$ in $L^1([a,b];U),$ one has that
$$
\label{xisunif}
(\xi_k^\tau),\, (\xi_k^{\tau,*})\  \text{converge uniformly to }\xi := \xi[\xib,u,v],
$$
with $\xib := \varphi(\xb,u(a)).$
Observe that, from previous equation and \eqref{uktau}-\eqref{tildeuktau}, one obtains
\benl
|x(\tau)-x^*(\tau)|=\lim |x_k^\tau(\tau)-x_k^{\tau,*} (\tau)|= \lim |\varphi \circ (\xi_k^\tau,-u^\tau_k)(\tau) - \varphi\circ (\xi_k^{\tau,*},-u^{\tau,*}_k)(\tau)| = 0.
\eenl
This contradicts our assumption, so the uniqueness of the limit solution of \eqref{E}-\eqref{IC} is proved. This ends the proof of Theorem  \ref{Existence}.

Finally, since $x$ was defined in  \eqref{xdef} via the coordinates' transformation, by  the uniqueness of the limit solution we also obtain the proof of Theorem \ref{TRep}.
\epf

{\em Proof of Theorem \ref{cont-depL1}.}
We shall start by proving item (ii). For this, let $r>0$ be arbitrary and consider the compact set $\tilde K'\subset \cR^{n+m}$ provided by the item (ii) of Theorem \ref{ContAC-simple}. Then, for each $|\xb|\leq r$ and $(u,v)\in \L^1([a,b];U)\times L^1([a,b];V),$  the corresponding limit solution $x[\xb,u,v]$ has values inside the compact set $K'$ defined as the projection of $\phi^{-1}(\tilde K')$ in the first $n$ components. This proves item (ii).

In order to show part (i), fix $\xb \in \cR^n$ and $u\in \L^1([a,b];U).$ Set $\xib := \varphi(\xb,u(a)),$ and notice that  $v \mapsto \xi[\xib,u,v]$ is continuous from $L^1([a,b];V)$ to $AC([a,b];\cR^n),$ where the latter space is endowed with the $\C^0-$topology. Furthermore, for any $v \in L^1([a,b];V),$ the trajectories $\xi[\xib,u,v]$ remain inside a compact set $\hat K.$
Let $v_k$ converge to $v$ in $L^1,$ then $\xi[\xib,u,v_k] \to \xi[\xib,u,v]$ uniformly and, since $\varphi$ is uniformly continuous on $\hat K \times (-U),$ one gets that
$$
v \mapsto x[\xb,u,v] = \varphi\circ (\xi[\xib,u,v],-u),
$$
is continuous from $L^1([a,b];V)$ to $L^\infty([a,b];\cR^n),$ from which  item (i) follows.

We  now prove item (iii). For an arbitrary $r>0,$  let $\tilde r>0$ be such that the points
$(\xib_1,u_1(a)):=\phi(\xb_1,u_1(a))$ and $(\xib_2, u_2(a)):=\phi(\xb_2, u_2(a))$ remain inside $\{|(\xib,\zetab)|\leq \tilde r\}$ whenever $|\xb_1|,|\xb_2| \leq r$ and $u_1,u_2 \in \L^1([a,b];U).$
 Take $\tau \in [a,b],$
 $|\xb_1|,|\xb_2| \leq r,$ $u_1, u_2 \in \L^1([a,b];U)$ and $v\in L^1([a,b];V).$ Let $x_1:=x[\xb_1,u_1,v]$ and $ x_2:=x[\xb_2, u_2,v]$ be the corresponding limit solutions of \eqref{E}.

 Let $(u_{1,k}^\tau),$ $(u_{2,k}^\tau)$ be sequences in $AC([a,b];U)$ verifying
 $$
 |u_{1,k}^\tau(a)-u_1(a)|+|u_{1,k}^\tau(\tau)-u_1(\tau)|+\|u_{1,k}^\tau -u_1 \|_1 \to 0,
 $$
 and
  $$
 | u_{2,k}^\tau(a)- u_2(a)|+| u_{2,k}^\tau(\tau)- u_2(\tau)|+\| u_{2,k}^\tau - u_2\|_1 \to 0,
 $$
 Consider the functions $\xi_{1,k}^\tau:=\xi[\xib_{1,k}^\tau,u_{1,k}^\tau,v],$ $\xi_{2,k}^\tau:=\xi[\xib_{2,k}^\tau, u_{2,k}^\tau,v],$ where $\xib_{1,k}^\tau:=\varphi(\xb_1,u_{1,k}^\tau(a)),$ $\xib_{2,k}^\tau:=\varphi(\xb_2, u_{2,k}^\tau(a)).$
 Then $\xi_{1,k}^\tau$ converge uniformly to $\xi_1:=\xi[\xib_1,u_1,v],$ where $\xib_1:=\varphi(\xb_1,u_1(a)),$ and $ \xib_{2,k}^\tau$ converge uniformly to $\xi_2:=\xi[ \xib_2,u_2,v],$ where $\xib_2:=\varphi(\xb_2, u_2(a)).$
 Let $\tilde M>1$ and $\tilde K'$
 be the constant and the compact set provided by Theorem \ref{ContAC-simple} for $\tilde r,$ respectively.
 Then one has
 \benl
 \begin{split}
& |\xi_{1,k}^\tau(\tau)-\xi_{2,k}^\tau(\tau)|+\|\xi_{1,k}^\tau-\xi_{2,k}^\tau\|_1 \leq \\
 &(\tilde M-1)\big[ |\xib_1-\xib_2| +  |u_{1,k}^\tau(a)- u_{2,k}^\tau(a)|+|u_{1,k}^\tau(\tau)- u_{2,k}^\tau(\tau)|+\|u_{1,k}^\tau - u_{2,k}^\tau\|_1 \big].
 \end{split}
 \eenl
 Letting $k$ go to infinity in the previous inequality, one obtains
\bel{xitildexi}
 \begin{split}
 &|\xi_1(\tau)-\xi_2(\tau)|+\|\xi_1- \xi_2\|_1\leq  \\
 &(\tilde M-1)\big[ |\xib_1-\xib_2| +  |u_1(a)- u_2(a)|+|u_1(\tau)- u_2(\tau)|+\|u_1 - u_2\|_1 \big].
 \end{split}
 \ee
 Thus, for $L$ a Lipschitz constant of $\varphi$ on $\tilde K',$ one gets
 \benl
 \begin{split}
 |x_1&(\tau)- x_2(\tau)|+\|x_1- x_2\|_1=\\
 &|\varphi\circ (\xi_1,-u_1)(\tau)-\varphi \circ ( \xi_2, u_2)(\tau)|
 +\|\varphi\circ (\xi_1,-u_1)-\varphi \circ ( \xi_2 , u_2 )\|_1 \leq \\
 &L\big[ |(\xi_1,u_1)(\tau)-( \xi_2 , u_2)(\tau)|+\|(\xi_1,u_1)-( \xi_2 , u_2)\|_1 \big] \leq\\
  &L(\tilde M-1)\big[|\xib_1-\xib_2|+|u_1(a)- u_2 (a)|+|u_1(\tau)- u_2(\tau)|+\|u_1- u_2\|_1\big] + \\
 &L\big[ |u_1(\tau)- u_2(\tau)|+\|u_1- u_2\|_1\big] \leq \\
  & M \big[|\xb_1- \xb_2|+|u_1(a)- u_2(a)|+|u_1(\tau)- u_2(\tau)|+\|u_1- u_2\|_1\big],
 \end{split}
 \eenl
 for $M>0$ depending on $L$ and $\tilde M,$ where the second inequality follows from \eqref{xitildexi}. Thus, the desired estimate holds true, and this completes the proof.
 \epf

{\em Proof of Theorem \ref{negl-jumps}.}
Set $\xib=\varphi(\xb,u(a)) = \varphi(\xb,\uh(a))$
and observe that $\xi[\xib,u,v](t)=\xi[\xib,\uh,v](t),$ for all $t\in [a,b].$ Let $\xi$ denote the latter function.
In view of Theorem \ref{TRep}, one has
\begin{gather*}
\xi(t) = \varphi \big( x[\xb,u,v](t),u(t) \big)=x[\xb,u,v](t) e^{- u_\alpha(t) g_\alpha}, \\
\xi(t) = \varphi \big( x[\xb,\uh,v](t),u(t) \big)=x[\xb,\uh,v](t)e^{- \uh_\alpha(t) g_\alpha}.
\end{gather*}
Therefore, combining the last two equations and due to the commutativity of $g_\alpha,$ the relation \eqref{jump} follows.
 \epf

\section{Concluding remarks}

A notion of  everywhere  defined solution  for the control Cauchy problem on $[a,b]$  $$ \dot{x} = {f}(t,x,u,v) +\sum_{\alpha=1}^m g_\alpha(x,u) \dot u_\alpha, \qquad t\in [a,b] \quad x(a)=\bar x,$$
has been provided,  under a commutativity hypothesis on the fields $g_\alpha$. In particular, we have proved results of existence, uniqueness and continuous dependence  on the data,  besides investigating  the effects of $u$'s changes on null sets. This concept of solution, which relies on an {\it extension by density} of the classical notion,  turns out to verify consistency requirements. We point out that, by defining the output at every $t\in [a,b],$  we have departed from a topological picture based on normed spaces, instead framing the limiting processes in spaces endowed with family of seminorms (this choice is concretely  represented by the fact  that approximating sequences in the solution's definition  depend on $\tau$, for every $\tau\in [a,b]$).

 The paper is  motivated by both applications (see the Introduction and Example \ref{minimum}) and the concern of constructing a suitable framework for further theoretical issues, like the study of  the corresponding {\it adjoint equations}, a likely crucial object in the investigation of necessary conditions for minima.

We think that a generalization of the notion of limit solution  to the  noncommutative case (in particular, an extension that will agree with former concepts of solutions) might  represent a natural direction for  further investigations. 

\section{Funding}

This work was partially supported by the European Union under the 7th Framework Programme FP7-PEOPLE-2010-ITN -  Grant agreement number 264735-SADCO, and the Fondazione CaRiPaRo Project 
``Nonlinear Partial Differential Equations: models, analysis, and 
control-theoretic problems".

\vspace*{6pt}


\begin{thebibliography}{99}


\bibitem{AruKarPer11}
\textsc{A.~Arutyunov, D.~Karamzin \& F.L. Pereira} (2011)
\newblock On a generalization of the impulsive control concept: controlling
  system jumps.
\newblock {\em Discrete Contin. Dyn. Syst.}, \textbf{29} (2), 403--415.


\bibitem{AruKarPer12}
\textsc{A.V. Arutyunov, D.Yu. Karamzin \& F.~Pereira} (2012)
 Pontryagin's maximum principle for constrained impulsive control
  problems, {\em Nonlinear Anal.,} \textbf{75} (3), 1045--1057.
  

\bibitem{BhaTiw09}
\textsc{S. Bhat \& P.K. Tiwari} (2009)   Controllability of spacecraft attitude using control moment
  gyroscopes,   {\em IEEE Transactions on Automatic Control,}   \textbf{54} (3), 585--590.

\bibitem{BrePicBook}
\textsc{A.~Bressan \& B.~Piccoli} (2007)
\newblock {\em Introduction to the mathematical theory of control}, vol.~2 of
  {American Institute of Mathematical Sciences (AIMS): AIMS Series on Applied Mathematics,}
\newblock Springfield, MO.


\bibitem{BreRam88}
\textsc{A.~Bressan \& F.~Rampazzo} (1988)
\newblock On differential systems with vector-valued impulsive controls,
\newblock {\em Boll. Un. Mat. Ital. B (7)}, \textbf{2} (3), 641--656.

\bibitem{BreRam91}
\textsc{A.~Bressan \& F.~Rampazzo} (1991)
\newblock Impulsive control systems with commutative vector fields,
\newblock {\em J. Optim. Theory Appl.}, \textbf{71} (1), 67--83.

\bibitem{BreAldo89}
\textsc{Aldo Bressan} (1989)
\newblock Hyper-impulsive motions and controllizable coordinates for
  {L}agrangian systems,
\newblock {\em Atti Accad. Naz. Lincei Mem. Cl. Sci. Fis. Mat. Natur. Sez. Ia
  (8)}, \textbf{19} (7), 195--246 (1991).

\bibitem{DalRam91}
\textsc{G.~Dal~Maso \& F.~Rampazzo} (1991)
\newblock On systems of ordinary differential equations with measures as
  controls,
\newblock {\em Differential Integral Equations}, \textbf{4} (4), 739--765.

\bibitem{Dyk94}
\textsc{V.A. Dykhta} (1994)
\newblock The variational maximum principle and quadratic conditions for the
  optimality of impulse and singular processes,
\newblock {\em Sibirsk. Mat. Zh.}, \textbf{35}  (1), 70--82, ii.


\bibitem{Haj85}
\textsc{O.~H\'{a}jec} (1985)
\newblock Book review,
\newblock {\em Bull. Amer. Math. Soc.,} \textbf{12} (2), 272--279.


\bibitem{Kar06}
\textsc{D.Yu. Karamzin} (2006)
\newblock Necessary conditions of the minimum in an impulse optimal control problem,
\newblock {\em J. Math. Sci.,} \textbf{139} (6), 7087--7150.

\bibitem{Lang}
\textsc{S.~Lang} (1995)
\newblock {\em Differential and {R}iemannian manifolds}, vol. 160 of {\em
  Graduate Texts in Mathematics},
\newblock Springer-Verlag, New York, third edition.

\bibitem{Ox80}
\textsc{J.~Oxtoby} (1980)
\newblock {\em Measure and Cathegory}, vol.~2 of {\em Graduate Texts in
  Mathematics}.
\newblock Springer-Verlag, New York, second edition.

\bibitem{Ram99}
\textsc{F.~Rampazzo} (1999)
\newblock Lie brackets and impulsive controls: an unavoidable connection,
\newblock In {\em Differential geometry and control ({B}oulder, {CO}, 1997)},
  volume~64 of {\em Proc. Sympos. Pure Math.}, pp. 279--296. Amer. Math.
  Soc., Providence, RI.

\bibitem{Ris65}
\textsc{R.W. Rishel} (1965)
\newblock An extended {P}ontryagin principle for control systems whose control
  laws contain measures.
\newblock {\em J. Soc. Indust. Appl. Math. Ser. A Control}, \textbf{3}, 191--205.

\bibitem{Sar91}
\textsc{A.V. Sarychev} (1991)
\newblock Nonlinear systems with impulsive and generalized function controls.
\newblock In {\em Nonlinear synthesis ({S}opron, 1989)}, vol.~9 of {\em
  Progr. Systems Control Theory}, pp. 244--257. Birkh\"auser Boston, Boston,
  MA.

\bibitem{SilVin96}
\textsc{G.N. Silva \& R.B. Vinter} (1996)
\newblock Measure driven differential inclusions.
\newblock {\em J. Math. Anal. Appl.}, \textbf{202} (3), 727--746.


\end{thebibliography}


\end{document}